\documentclass{amsart}
\usepackage{latexsym,amsxtra,amscd,ifthen}
\usepackage{amsfonts}
\usepackage{verbatim}
\usepackage{amsmath}
\usepackage{amsthm}
\usepackage{amssymb}

\numberwithin{equation}{section}

\usepackage{amssymb,amsmath}
\usepackage{xcolor}

\newtheorem{Theorem}{Theorem}

\newtheorem{Proposition}[Theorem]{Proposition}

\newtheorem{Conjecture}[Theorem]{Conjecture}

\newcommand{\kk}{{\Bbbk}} 

\newcommand{\Complex}{{\mathbb C}}

\def\ep{\hfill{\vbox to 7pt{\hbox to 7pt{\vrule height 7pt width 7pt}}}}

\DeclareMathOperator{\AdOper}{Ad}

\newcommand{\Identity}{Id}

\begin{document}
	
	\title{Multilinear nilalgebras and the Jacobian theorem}

	\author{Dmitri Piontkovski}
	
	\address{20 Myasnitskaya ul., National Research University Higher School of Economics, Moscow 101000, Russia\\
ORCID: https://orcid.org/0000-0002-9853-1891
	}
	
	\email{dpiontkovski@hse.ru}
	
	

	\date{\today}
	
  \maketitle   	
 
\begin{abstract}
If a symmetric multilinear algebra is weakly nil, then it is Engel.  This result may be regarded as an infinite-dimensional analogue of the well-known Jacobian theorem, which states that if a polynomial mapping has a polynomial inverse, then its Jacobian matrix is invertible. This refines a theorem of Gerstenhaber and partially answers a question posed by Dotsenko. 
\end{abstract}

 \label{sec:Jacobian}

 There are various versions of nil algebras in non-associative and multilinear settings. These notions are closely connected with problems concerning polynomial automorphisms.  

 Let $A$ be a vector space and $\mu: A^m \to A$ be a $m$-linear map, so that the pair $(A,\mu)$ forms a (multilinear) algebra (in the sense of general algebra).
The algebra $A$ is called {\em $n$-Engel} (or simply Engel) if for each element $x\in A$ the linear operator of multiplication by $x$ is $n$-nilpotent, that is, the $n$-th Engel identity 
$$
   \AdOper_x^n(y) \equiv 0 
$$ 
holds in $A$, where $\AdOper^1_x (y)  = \AdOper_x (y) = \mu (x, \dots, x, y), 
 \AdOper^{k+1}_x (y) = \mu (x, \dots, x, \AdOper^{k}_x (y))$. 

Gerstenhaber  connected the Engel property with another kind of nilpotence. We call $A$ {\em Gerstenhaber nil} (of nilindex $n$) if for each element $x\in A$, all multilinear multiple compositions of $\mu$ with at least  $n$ arguments each of them equal to $x$, are zero.  In other words, $A$ is Gerstenhaber nil if each subalgebra generated by a single element is nilpotent of bounded degree. 

From now on, assume that   the  ground field $\kk$ has zero characteristic and the algebra $A$ is commutative in the sense that the operation $\mu$ is symmetric, i.~e., for every $\sigma \in S_m$ we have
 	$\mu(x_1, \dots , x_m ) = \mu(x_{\sigma 1}, \dots , x_{\sigma m})$. 
The following theorem was proved by Gerstenhaber~\cite{Gerst1960}  for binary algebras and generalized  to the general case $m\ge 2$ by Umirbaev~\cite[Lemma 4.1]{Umirbaev}. 

\begin{Theorem}[Gerstenhaber, Umirbaev]
Each Gerstenhaber nil algebra is Engel. 
\end{Theorem}
In the binary case Gerstenhaber proved the following  estimate: if $A$ is Gerstenhaber $p$-nil, then it is $n$-Engel for  $n = 2p-3$~\cite[Theorem~1]{Gerst1960}.  For a generalization of this estimate to the case of $m$-linear algebras in a more general setting, see Theorem~\ref{th:main} below.

Another version of nilpotence appears in an equivalent formulation of the Jacobian problem.
 Let $T^{mult}_q (x_1, \dots , x_q)$ be the sum of all multilinear multiple compositions of $\mu$ with $q$ arguments, and let $T_q(x) = \frac{1}{q!}T_q^{mult}(x,x, \dots, x)$. 
In other words, 
$$
T_1(x) = x \mbox{ and }T_q (x) = \sum_{i_1 +\dots + i_m = q}\mu(T_{i_1}, \dots, T_{i_m}) 
\mbox{ for }q>1. 
$$ 

Let us call an algebra $A$ {\em Yagzhev nil}  if the identities $T_q(x) \equiv 0$ hold in $A$ for all $q\ge p$. The minimal such $p$ is called the (Yagzhev) {\em  nilindex} of $A$. Such algebras are called Yagzhev, or weakly nilpotent  in~\cite{Belov_etc}, and weakly nil in~\cite{dots}.  Note that Yagzhev nil algebras need not be Gerstenhaber nil.  Indeed, 
Gorni and Zampieri~\cite{gz} have shown that there exists a 4-dimensional 3-linear algebra which is Yagzhev nil and Engel but not Gerstenhaber nil (their example is induced by a polynomial automorphism due to van den Essen).  However, in the important case of binary algebras we do not know whether these two nil properties are equivalent? At least, a straightforward calculation with linearized identities shows that binary Yagzhev nil  algebras of nilindex $4$ are Gerstenhaber nil of nilindex $6$.

Yagzhev nil algebras  appear in his reformulation of the famous Jacobian conjecture, see~\cite{Belov_etc, Umirbaev} and references therein.  Recall that the conjecture states that a complex polynomial endomorphism with constant Jacobian determinant has a polynomial inverse. 
 It is equivalent to the following statement (\cite{Yagzhev}; see also \cite{bass}):
 
 \begin{Conjecture}[Jacobian conjecture for homogeneous mapping]
 	\label{conj:Jacobian_hom}
 	Suppose that $F: \Complex^n \to \Complex^n$ is a polynomial automorphism of the form 
 	$F = \Identity - H$, where $H$ is a homogeneous automorphism of degree $m\ge 2$.
 	If the determinant $j_F$ of the Jacobian matrix $J_F = (\partial F(X)/\partial X)$ of the map $F$ is equal to the constant 1, then the map $F$ has a polynomial inverse. 
 \end{Conjecture}
 
 The conjecture is known to be true for $m=2$~\cite{Wang}, and it is sufficient to prove it for $m=3$~(\cite{Yagzhev,bass}). 
 
 It follows from the results by Yagzhev~\cite{Yagzhev2000} (based on an earlier work by Dru\. zkowski and Rusek~\cite{dr}) 
that the Jacobian conjecture for homogeneous mapping of degree $m$ for polynomials of $n$ variables  is equivalent to the following:
 
 \begin{Conjecture}[Jacobian conjecture in Yagzhev form]
 	\label{conj:Yagzhev}
 	Suppose that $A$ is a complex $n$-dimensional  $m$-linear algebra. If $A$ is Engel, 
then it is Yagzhev nil.  
 \end{Conjecture}
 In the notation of Conjecture~\ref{conj:Jacobian_hom}, here $\mu$ is the complete linearization of $H$ so that $ H(X) =\mu(X, \dots, X)$.  


 The easy reverse implication of the Jacobian conjecture (stating that the Jacobian determinant of a polynomial automorphism is a nonzero constant) is  well known as the  Jacobian theorem. Being reformulated in the Yagzhev terms, it is a strengthening of the Gerstenhaber--Umirbaev theorem for finite-dimensional algebras.

 \begin{Proposition}[Jacobian theorem]
 	\label{prop:Easy_part_of_JC}
 	Suppose that a  complex $m$-linear algebra $A$ is finite-dimensional.
If $A$ is Yagzhev nil, then it is  Engel. 
 \end{Proposition}
 
The aim of this note is to prove the same result without the assumption that $A$ is finite-dimensional. We regard this as an infinite-dimensional version of the Jacobian theorem.

 \begin{Theorem}[Jacobian theorem for an infinite set of variables]
 	\label{th:main}
	Each $m$-linear Yagzhev nil algebra $A$ is  Engel. 
If the Yagzhev nilindex of $A$ is $p$, 
then $A$ is $n$-Engel
for $n =p-1+\frac{p-2}{m-1} $. 
 \end{Theorem}
Note that the assumption that the algebra $A$ is Yagzhev nil of index $p$ in the second implication  can be replaced by the following: the identities $T_q(x) \equiv 0$ hold for all $p \le q \le m(p-1)+1$ (cf.~\cite[Prop.~3]{dots}). 

 Dotsenko has asked whether, in the binary case $m=2$, the identity $T_q(x) \equiv 0$  implies the $n$-Engel identity  for some $n$ \cite[Question 4]{dots}.
Theorem~\ref{th:main} gives a partial answer to this question. 
Indeed,  if the identities $T_q(x) \equiv 0$ hold for all $p \le q \le 2p-1$, the algebra is $n$--Engel for $n = 2p-3$.  Moreover, one can show that if the identities $T_4(x) \equiv T_5(x)  \equiv 0$ hold in a binary algebra $A$, then it is Gerstenhaber nil  of nilindex 6 (hence, Yagzhev nil of nilindex 4); therefore, it is 5-Engel (cf.~\cite[Question~3]{dots}).

\medskip

{\it Proof of Theorem~\ref{th:main}. }\ 

 We may (and will) assume that all relations of the algebra $A$ are consequences 
 of the  
 identities $T_q(x) \equiv 0$ for all $q\ge p$, so that $A$ is a free algebra of the variety $V$ of $m$-linear  algebras defined by these identities. Moreover, we will assume that the set of free  generators of $A$ consists of at least two elements. Then $A$ admits a natural grading such that the degrees of the free  generators are 1. 

 The operator $g:A\to A$ defined by $g(x) = x - \mu (x, \dots ,x)$ is invertible; 
 the  inverse is given by $\gamma(y) := g^{-1}(y) = \sum_{j\ge 1} T_j(y)$~\cite{dr} (where $T_j(y)=0$ for $j\ge p$).
 We obtain the identities
 $$
 \gamma (g(x)) = \sum_{j\ge 1} T_j(g(x)) \equiv  x
 \mbox{ and }  g(\gamma(y)) = \gamma(y)-\mu (\gamma(y), \dots ,\gamma(y)) \equiv y.
 $$
Replacing $x$ by $x+z$, $y$ by $y+t$, and collecting all terms which are linear in $z$ and $t$, we obtain partial linearizations of the above identities
 $$
 \mathrm{d}\gamma (g(x), \mathrm{d}g(x,z))  \equiv  z\mbox{ and }
 \mathrm{d}g (\gamma(y), \mathrm{d}\gamma(y,t))  \equiv  t,
 $$
 where 
$$\mathrm{d}\gamma(y,t) = \sum_{n= 1}^{p-1} \frac{1}{(j-1)!}T_j^{mult}(y,\dots, y,t)
\mbox{ and }
\mathrm{d}g(x,z) = z - m \AdOper_x(z)
$$ are the corresponding partial linearizations of $\gamma$ and $g$. 
 So, the linear operator $\mathrm{d}g_x: z \mapsto \mathrm{d}g(x,z)$
on $A$
 (the ``Jacobian'') is invertible; the inverse is given by 
 the operator $t \mapsto \mathrm{d}\gamma(g(x),t)$.
 
 On the other hand, 
the linear operator $\mathrm{d}g_x$ can be extended to the 
completion $\widehat A$ of the graded algebra $A$ 
by the same formula $\widehat{\mathrm{d}g_x}: z\mapsto z -m \mu(x, \dots, x,z)$.
It admits the inverse  
 $\widehat{\mathrm{d}g_x}^{-1}: t \mapsto t+ \sum_{i\ge 1} m^i \AdOper_x^i(t)$. Since 
 the restriction of the operator $\widehat{\mathrm{d}g_x}$ to the subset $A\subset \widehat A$ is a bijection $\mathrm{d}g_x:A\to A$, we have $\widehat{\mathrm{d}g_x}^{-1}(A)\subset A$. 
As $A$ is naturally graded, only finitely many
homogeneous components $m^i \AdOper_x^i(t)$ of $\widehat{\mathrm{d}g_x}^{-1}(t)$
 are nonzero.  That is, $\AdOper_x^i(t)$ for $i>>0$.

Assume now that that $x$ and  $t$ are different free generators of $A$.  
Let $N$ be the smallest number such that  the equalities  
$\AdOper_x^i(t) =0$ hold for all $i\ge N$. Since $A$ is a free algebra of the variety $V$,
these equaltites are identities of  $V$. In particular, $A$ is $N$-Engel. 

To evaluate  $N$, consider the
equality
$$
\widehat{\mathrm{d}g_x}^{-1}(t) = {\mathrm{d}g_x}^{-1}(t)
$$
for the free generators $x,t$ of $A$. 
Since the algebra is graded, this equality implies the equality of corresponding homogeneous components. On the left-hand side  
$\widehat{\mathrm{d}g_x}^{-1}(t) =  t+ \sum_{i\ge 1} m^i \AdOper_x^i(t)$,
the highest nonzero component is the one with $i=N-1$; its degree is $(m-1)(N-1)+1$.
On the right-hand side  
$${\mathrm{d}g_x}^{-1}(t) = \mathrm{d}\gamma(g(x),t) = 
\sum_{j= 1}^{p-1} \frac{1}{(j-1)!}T_j^{mult}(g(x),\dots, g(x),t),
$$
the degrees of the nonzero homogeneous components
do not exceed the number $m(p-2) +1$. Therefore, we obtain the inequality
$$
(m-1)(N-1)+1 \le m(p-2) +1.
$$
Thus,   $N\le m\frac{p-2}{m-1} +1 = p-1+\frac{p-2}{m-1}$.  
In other words, $A$ is $n$-Engel for all $n\ge p-1+\frac{p-2}{m-1}$. 
 \ep\\
 
The converse of Theorem~\ref{th:main}  (whether an arbitrary Engel algebra is Yagzhev nil?) is a challenging open problem which generalizes the Jacobian conjecture to the case of infinite number of variables.  For binary algebras, it is stated in~\cite{Belov_etc} as the Generalized Jacobian conjecture for quadratic mappings.  
The conjecture is known to hold for several classes of algebras, including power-associative algebras, subalgebras of Bernstein algebras, algebras satisfying the identity$(x^2)^2\equiv 0$~\cite[Sec.~8]{pz}, and 3-Engel algebras~\cite{singer}. 
Note that in all these cases the algebras turn out to be not only Yagzhev nil but also Gerstenhaber nil.  

\medskip
 {\large \bf Acknowledgments}

I am grateful to Vladimir Dotsenko and Fouad Zitan for stimulating discussions.  My special thanks to the anonymous referee for detailed comments and suggestions. 
 This work was supported by the Russian Science Foundation grant 24-21-00341,  https://rscf.ru/project/24-21-00341/.

%

\medskip
\noindent{\large\bfseries Addition.}
\emph{This addition was made after the paper had been accepted
for publication. It appears only in the arXiv version and is
not included in the journal version.}

The Jacobian conjecture has recently been disproved by
Alp\"oge with AI assistance. The generalized Jacobian conjecture
for quadratic mappings, however, remains open.

We also present here an example (carried out with AI
assistance) which shows that for binary algebras, the Engel identity does not imply Gerstenhaber
nilpotence (for ternary algebras this was discussed above). The example is Umirbaev's modification of the
automorphism associated with Suttles's algebra; the relevant
nonnilpotence phenomenon was already noted
in~\cite[Sec.~5]{Umirbaev}. Namely, consider the polynomial
automorphism
\[
  F(X)=X-H(X),
  \qquad
  H(a,b,c,d,e)=(bd,\ ae-cd,\ be,\ -e^2,\ 0).
\]
Suppose that the characteristic of the ground field is not equal to 2. Then  the corresponding
five-dimensional polarization algebra $B$, with multiplication
\[
  X\cdot Y
  =\frac12\bigl(H(X+Y)-H(X)-H(Y)\bigr),
\]
has Engel index exactly $5$ and is Yagzhev nil of index $6$,
but is not Gerstenhaber nil. 

To prove the last assertion, let $e_1,\ldots,e_5$ denote the standard coordinate
basis of $B$, and put $x=e_2+e_5$. The algebra $B$ is generated
by $x$. Moreover, writing $\AdOper_y(z)=y\cdot z$ and
$x^3=(x\cdot x)\cdot x$, we obtain
\[
  x^2=e_3-e_4,
  \qquad
  x^3=-\frac12e_1,
  \qquad
  \AdOper_{x^2}\AdOper_x(e_1)=-\frac14e_1.
\]
Consequently,
\[
  \bigl(\AdOper_{x^2}\AdOper_x\bigr)^n(x^3)
  =\frac{(-1)^{n+1}}{2^{2n+1}}\,e_1
  \ne0
  \qquad(n\ge0).
\]
These are one-variable monomials of degree $3n+3$, so their
nonvanishing proves that $B$ is not Gerstenhaber nil. 

As a corollary, the free monogenic algebra in the variety of commutative binary algebras defined by the 5-Engel identity is infinite-dimensional.

 \end{document}